\documentclass[11pt]{article}
\usepackage{epsfig,amsmath,cite}
\usepackage{amssymb}
\usepackage{color}
\usepackage{amsfonts}

\def\be{\begin{equation}}
\def\ee{\end{equation}}
\def\bea{\begin{eqnarray}}
\def\eea{\end{eqnarray}}
\def\bes{\begin{eqnarray*}}
\def\ees{\end{eqnarray*}}

\def\nn{\nonumber}
\def\lb{\label}
\def\bs{\setminus}

\def\R{{\bf R}}
\def\C{{\bf C}}
\def\Z{{\bf Z}}
\def\K{{\bf K}}
\def\N{{\bf N}}
\def\U{{\bf U}}

\def\Q{{\bf Q}}

\def\HP{{\bf HP}}
\def\CaP{{\bf CaP}}
\def\CP{{\bf CP}}

\def\aa{{\alpha}}
\def\bb{{\beta}}
\def\ga{{\gamma}}

\def\th{{\theta}}
\def\om{{\omega}}
\def\Om{{\Omega}}
\def\ep{{\epsilon}}
\def\lm{{\lambda}}
\def\Lm{{\Lambda}}

\def\sg{{\sigma}}
\def\dm{{\diamond}}

\def\vf{{\varphi}}

\def\K{{\cal K}}
\def\P{{\cal P}}

\def\M{{\cal M}}
\def\Nn{{\cal N}}

\def\rank{{\rm rank}}
\def\Sp{{\rm Sp}}

\def\mod{{\rm mod}}
\def\ol{\overline}

\def\hb{\vrule height0.18cm width0.14cm $\,$}

\title{Multiplicity of closed geodesics on bumpy Finsler manifolds with elliptic closed geodesics}

\author{Huagui Duan\thanks{Partially supported by National Key R\&D Program of China (Grant No. 2020YFA0713300),
NNSFC (Nos. 12271268, 11671215 and 11790271) and the Fundamental Research Funds for the Central Universities. E-mail: duanhg@nankai.edu.cn.}, \qquad
Dong Xie\thanks{Partially supported by NNSFC (Nos. 12271268 and 11790271). E-mail: 2120170062@mail.nankai.edu.cn} \\ \\
School of Mathematical Sciences and LPMC, Nankai University,\\ Tianjin 300071, P. R. China\\}

\begin{document}
\maketitle

\begin{abstract}
{\it Let $M$ be a compact simply connected manifold satisfying $H^*(M;\Q)\cong T_{d,n+1}(x)$ for
integers $d\ge 2$ and $n\ge 1$. If all prime closed geodesics on $(M,F)$ with an irreversible bumpy Finsler metric $F$ are elliptic, then
either there exist exactly $\frac{dn(n+1)}{2}$ (when $d\ge 2$ is even) or $(d+1)$ (when $d\ge 3$ is odd) distinct closed
geodesics, or there exist infinitely many distinct closed geodesics. }
\end{abstract}



\renewcommand{\theequation}{\thesection.\arabic{equation}}
\renewcommand{\thefigure}{\thesection.\arabic{figure}}

\setcounter{equation}{0}
\section{Introduction and main result}

Recall that a closed curve on a Finsler manifold is a closed geodesic if it is locally the shortest path connecting any two nearby points on this curve (cf. \cite{She01}). As usual, on any Finsler manifold $(M, F)$, a closed geodesic $c:S^1=\R/\Z\to M$ is {\it prime}
if it  is not a multiple covering (or, iteration) of any other closed geodesic. Here the $m$-th iteration $c^m$ of $c$ is defined
by $c^m(t)=c(mt),\ \forall\,m\ge 1$.  The inverse curve $c^{-1}$ of $c$ is defined by $c^{-1}(t)=c(1-t)$ for $t\in \R$. In the Riemannian case, two closed geodesics
$c$ and $d$ are called {\it geometrically distinct} if $c(S^1)\neq d(S^1)$. However, unlike the Riemannian case, the inverse curve $c^{-1}$ of a closed geodesic $c$ on an irreversible Finsler manifold $(M,F)$ need not be a geodesic. And we call two prime closed geodesics $c$ and $d$ on $(M,F)$ {\it distinct} if there is no $\th\in (0,1)$ such that $c(t)=d(t+\th)$ for all $t\in\R$.

For a closed geodesic $c$ on an $(n+1)$-dimensional manifold $M$, denote by $P_c$ the linearized Poincar\'{e} map of $c$, which is a symplectic matrix, i.e., $P_c\in\Sp(2n)$. We define the {\it elliptic height } $e(P_c)$ of $P_c$ to be the total algebraic multiplicity of all eigenvalues of $P_c$ on the unit circle $\U=\{z\in\C|\; |z|=1\}$ in the complex plane $\C$. Since $P_c$ is symplectic, $e(P_c)$ is even and $0\le e(P_c)\le 2n$.
A closed geodesic $c$ is called {\it elliptic} if $e(P_c)=2n$, i.e., all the eigenvalues of $P_c$ locate on $\U$; {\it irrationally elliptic} if, in the homotopy component $\Om^0(P_c)$ of $P_c$ (cf. Section 2 below for the definition), $P_c$ can be connected to the $\dm$-product of $n$ rotation matrices $R(\th_i)$ with $\th_i$ being irrational multiple of $\pi$ for $1\le i\le n$; {\it hyperbolic} if $e(P_c)=0$, i.e., all the
eigenvalues of $P_c$ locate away from $\U$; {\it non-degenerate} if $1$ is not an eigenvalue of $P_c$. A Finsler metric $F$
is called {\it bumpy} if all the closed geodesics on $(M,F)$ are non-degenerate.

The closed geodesic problem is a classical one in Riemannian geometry and dynamical system. A famous conjecture claims the existence of infinitely
many distinct closed geodesics on every compact Riemannian manifold $(M,g)$. This conjecture has been proved for many cases. For example,
Gromoll and Meyer proved the conjecture in \cite{GrM69} for any compact $M$ provided the Betti number sequence $\{b_p(\Lm M)\}_{p\in\N}$ of the free
loop space $\Lm M$ of $M$ is unbounded. Later, for any compact simply connected manifold $M$, Vigu\'e-Poirrier and Sullivan in \cite{ViS76} further proved that the Betti number sequence is bounded if and only if $M$ satisfies \bea  H^*(M;\Q)\cong T_{d,n+1}(x)=\Q[x]/(x^{n+1}=0)   \lb{1.1}\eea
with a generator $x$ of degree $d\ge 2$ and height $n+1\ge 2$, where $\dim M=dn$. Note that when $d$ is odd, then $x^2=0$ and $n=1$, and then $M$ is rationally homotopic to $S^d$ (cf. Remark 2.5 of \cite{Rad89} and \cite{Hin84}). Especially, these manifolds include the sphere $S^d$ (with $n=1$), complex projective space $\CP^n$ (with $d=2$), quaternionic projective space $\HP^n$ (with $d=4$), and Cayley plane $\CaP^2$ (with $d=8$ and $n=2$) as special examples.

Among manifolds satisfying (\ref{1.1}), only for Riemannian $S^2$ the above conjecture was proved by \cite{Fra92} of Franks and \cite{Ban93} of Bangert. When considering Finsler metrics, the situation changes dramatically. It was quite surprising that Katok in \cite{Kat73} constructed some irreversible Finsler metrics on rank one symmetric spaces which possess finitely many ($\frac{dn(n+1)}{2}$ when
$d\ge 2$ is even, or $d+1$ when $d\ge 2$ is odd) distinct closed geodesics, and all these closed geodesics are irrationally elliptic. The geometry of Katok's metrics was further studied by Ziller in \cite{Zil82}. Based on Katok's examples, Anosov in \cite{Ano74} conjectured that the smallest number of distinct closed geodesics on any Finsler sphere $S^d$ should be $2[(d+1)/2]$. This conjecture was proved by Bangert and Long in \cite{BL10} for every Finsler sphere $(S^2,F)$.

It is natural to generalize this conjecture to
all compact simply connected Finsler manifolds satisfying (\ref{1.1}), i.e., the lower bound
of the number of distinct closed geodesics on such manifolds should be $\frac{dn(n+1)}{2}$ when
$d\ge 2$ is even, or $d+1$ when $d\ge 2$ is odd, respectively. This has been proved by \cite{DLW16a} under the bumpy metric and the nonnegative flag curvature conditions.

Recently, a great number of results about the multiplicity and
stability of closed geodesics on Finsler manifolds have appeared, for example, we refer
readers to some recent survey \cite{BuK21}, \cite{DLZ20}, \cite{Lon06}, \cite{Tai10}, and some research papers \cite{Dua15}, \cite{Dua16}, \cite{DL22}, \cite{DuL07}, \cite{DuL10}, \cite{DLW16b}, \cite{GG20}, \cite{GGM18}, \cite{GM21}, \cite{HR13}, \cite{HWZ98}, \cite{HWZ03}, \cite{Liu17}, \cite{LX17}, \cite{LLX18}, \cite{LD09}, \cite{Rad07}, \cite{Rad10}, \cite{Wan08},\cite{Wan12}, \cite{Wan13} and the references therein.

Note that Hingston in \cite{Hin84} proved the existence of infinitely many closed geodesics on spheres if all prime closed geodesics are hyperbolic. Hofer, Wysocki and Zehnder in \cite{HWZ98} and \cite{HWZ03} proved that there exist either two or infinitely many prime closed geodesics on every bumpy Finsler $S^2$ provided the stable and unstable manifolds of every hyperbolic closed geodesic intersect transversally. Furthermore, most recently this result has been proved to be true only under the the bumpy metric by \cite{CHP19} of Cristofaro-Gardiner, Hutchings and Pomerleano by using the embedded contact homology method. In addition, Long and Wang in \cite{LW08} proved that on every Finsler $S^2$, there exist either at least two irrationally elliptic closed geodesics or infinitely many closed geodesics.

Based on Katok's examples and the above well-known results, about the number of closed geodesics on compact simply-connected
manifolds with irreversible Finsler metrics, we suspect the following result holds.

\medskip

{\bf Conjecture A.} (cf. \cite{DLW16a}) {\it Let $M$ be a compact simply connected manifold satisfying (\ref{1.1}), and $F$ be any irreversible Finsler metric $F$ on $M$. Denote by $\Nn(M,F)$ the number of the distinct closed geodesics on $(M,F)$. Then there holds
\begin{itemize}
\item $\Nn(M,F) = \frac{dn(n+1)}{2}\;$ or $\;+\infty$, if $d$ is even;
\item $\Nn(M,F) = d+1\;$ or $\;+\infty$, if $d$ is odd.
\end{itemize}
Furthermore, if $\Nn(M,F)<+\infty$, then all these closed geodesics are irrationally elliptic.}

\medskip

Motivated by Katok's examples and the above Conjecture, we will establish the following result, which gives a partial answer to Conjecture A in the bumpy case.

\medskip

{\bf Theorem 1.1.} {\it Let $M$ be a compact simply connected manifold satisfying (\ref{1.1}). If all prime closed geodesics on $(M,F)$ with an irreversible bumpy Finsler metric $F$ are elliptic, then either $\Nn(M,F) \in \{\frac{dn(n+1)}{2},\,+\infty\}$ when $d$ is even, or $\Nn(M,F) \in \{ d+1,\;+\infty\}$ when $d$ is odd.}

\medskip

Note that when the manifold $(M,F)$ is the sphere $(S^d,F)$, Theorem 1.1 has been proved in \cite{DL16} under more stronger conditions. On the other hand, the methods in \cite{DL16} depend heavily on the monotonicity of iterated indices of closed geodesics, which has been guaranteed since the Morse index $i(c)\ge \dim M-1$ there for every closed geodesic $c$.

In this paper, note that $\dim M=dn$, but the Morse index can only be estimated to be greater than and equal to $(d-1)$ (see Claim 2 in Section 3), which will lead to the loss of monotonicity of  the iterated indices of closed geodesics. In order to overcome this difficulty, we will make use of some methods and arguments developed most recently by the first author and his co-workers, especially the enhanced common index jump methods established in \cite{DLW16a} and \cite{DLW16b}.

This paper is organized as follows. In Section 2, we review the critical point theory and index iteration theory for closed geodesics. In Section 3, we give the proof of Theorem 1.1. Next let $\N$, $\N_0$, $\Z$, $\Q$, $\R$, and $\C$ denote the sets of natural integers, non-negative integers, integers, rational numbers, real numbers, and complex numbers respectively. We use only singular homology modules with $\Q$-coefficients.
The following notations will be used in this paper.
\bea \left\{\begin{array}{ll}[a]=\max\{k\in\Z\,|\,k\le a\}, &
           E(a)=\min\{k\in\Z\,|\,k\ge a\} , \cr
    \varphi(a)=E(a)-[a], &\{a\}=a-[a]. \\ \end{array} \right. \lb{1.2}\eea
Especially, $\varphi(a)=0$ if $ a\in\Z\,$, and $\varphi(a)=1$ if $
a\notin\Z\,$.


\setcounter{equation}{0}
\section{Morse theory and index iteration theory of closed geodesics}
\subsection{Morse theory of closed geodesics}
Let $(M,F)$ be a compact Finsler manifold, the space
$\Lambda=\Lambda M$ of $H^1$-maps $\gamma:S^1\rightarrow M$ has a
natural structure of Hilbert manifold on which the
group $S^1=\R/\Z$ acts continuously by isometries (cf. \cite{Kli78}). This action is defined by
$(s\cdot\gamma)(t)=\gamma(t+s)$ for all $\gamma\in\Lm$ and $s,
t\in S^1$. For any $\gamma\in\Lambda$, the energy functional is
defined by
\be E(\gamma)=\frac{1}{2}\int_{S^1}F(\gamma(t),\dot{\gamma}(t))^2dt.
\lb{2.1}\ee
It is $C^{1,1}$ and invariant under the $S^1$-action. The
critical points of $E$ of positive energies are precisely the closed geodesics
$\gamma:S^1\to M$. The index form of the functional $E$ is well
defined along any closed geodesic $c$ on $M$, which we denote by
$E''(c)$. As usual, we denote by $i(c)$ and
$\nu(c)$ the Morse index and nullity of $E$ at $c$. In the
following, we denote by
\be \Lm^\kappa=\{d\in \Lm\;|\;E(d)\le\kappa\},\quad \Lm^{\kappa-}=\{d\in \Lm\;|\; E(d)<\kappa\},
  \quad \forall \kappa\ge 0. \nn\ee
For a closed geodesic $c$ we set $ \Lm(c)=\{\ga\in\Lm\mid E(\ga)<E(c)\}$.

For $m\in\N$ we denote the $m$-fold iteration map
$\phi_m:\Lambda\rightarrow\Lambda$ by $\phi_m(\ga)(t)=\ga(mt)$, for all
$\,\ga\in\Lm, t\in S^1$, as well as $\ga^m=\phi_m(\gamma)$. If $\gamma\in\Lambda$
is not constant then the multiplicity $m(\gamma)$ of $\gamma$ is the order of the
isotropy group $\{s\in S^1\mid s\cdot\gamma=\gamma\}$. For a closed geodesic $c$,
the mean index $\hat{i}(c)$ is defined as usual by
$\hat{i}(c)=\lim\limits_{m\to\infty}i(c^m)/m$. Using singular homology with rational
coefficients we consider the following critical $\Q$-module of a closed geodesic
$c\in\Lambda$:
\be \overline{C}_*(E,c)
   = H_*\left((\Lm(c)\cup S^1\cdot c)/S^1,\Lm(c)/S^1\right). \lb{2.2}\ee

{\bf Proposition 2.1.} (cf. Satz 6.11 of \cite{Rad92} ) {\it Let $c$ be a
prime closed geodesic on a bumpy Finsler manifold $(M,F)$. Then there holds}
$$ \overline{C}_q( E,c^m) = \left\{\begin{array}{ll}
     \Q, &\quad {\it if}\;\; i(c^m)-i(c)\in 2\Z\;\;{\it and}\;\;
                   q=i(c^m),\;  \cr
     0, &\quad {\it otherwise}. \\ \end{array}\right.  $$

{\bf Definition 2.2.} (cf. Definition 1.6 of \cite{Rad89}) {\it For a
closed geodesic $c$, let $\ga_c\in\{\pm\frac{1}{2},\pm1\}$ be the
invariant defined by $\ga_c>0$ if and only if $i(c)$ is even, and
$|\ga_c|=1$ if and only if $i(c^2)-i(c)$ is even. }

{\bf Theorem 2.3.} (cf. Theorem 3.1 of \cite{Rad89} and Satz 7.9 of \cite{Rad92}) {\it Let
$(M,F)$ be a compact simply connected bumpy Finsler manifold of $\dim M=dn$ satisfying
$\,H^{\ast}(M,\Q)=T_{d,n+1}(x)$. Denote by $\{c_k\}_{1\le k\le q}$ the prime closed geodesics on $(M,F)$
with positive mean indices. Then
\be \sum_{j=1}^q\frac{\ga_{c_j}}{\hat{i}(c_j)}=B(d,n)
=\left\{\begin{array}{cc}
     -\frac{n(n+1)d}{2d(n+1)-4}, &\quad d\;\;{\it is\;even},\cr
     \frac{d+1}{2d-2}, &\quad d\;\;{\it is\;odd}. \end{array}\right.   \lb{2.3}\ee
}

Let $(X,Y)$ be a space pair with $S^1$-action such that the Betti numbers $b_i=\dim H_i(X/S^1,Y/S^1;\Q)$ are finite for all $i\in \Z$. The following result gives the precise information about the Betti number sequence of $\Lm M$.

\medskip

{\bf Lemma 2.4.} (cf. Theorem 2.4 of \cite{Rad89} and Lemma 2.5 and Lemma 2.6 of \cite{DuL10})
{\it Let $M$ be a compact simply connected manifold satisfying $H^*(M;\Q)\cong T_{d,n+1}(x)$ for
integers $d\ge 2$ and $n\ge 1$.

(i) When $d$ is odd (which implies that $n=1$), i.e. $M$ is rationally homotopic to the sphere $S^d$, then the Betti numbers of the free loop space of $S^d$ are given by
\bea b_i(\Lm S^d)
&=& \rank H_i(\Lm S^d/S^1,\Lm^0 S^d/S^1;\Q)  \nn\\
&=& \left\{\begin{array}{ccc}
    2,&\quad {\it if}\quad i\in \K\equiv \{k(d-1)\,|\,2\le k\in\N\},  \cr
    1,&\quad {\it if}\quad i\in \{d-1+2k\,|\,k\in\N_0\}\bs\K,  \cr
    0 &\quad {\it otherwise}. \end{array}\right. \lb{2.4}\eea
For any integer $k\ge d-1$, there holds
\bea \sum_{i=0}^k b_i(\Lm S^d) = \left[\frac{k}{d-1}\right] + \left[\frac{k}{2}\right] - \frac{d-1}{2}.   \lb{2.5}\eea

 (ii) When $d$ is even, let $D=d(n+1)-2$ and
\bea &&\Om(d,n) = \{k\in 2\N-1\,|\, k_1D\le k-(d-1)=k_1D+k_2d\le k_1D+(n-1)d\;  \nn\\
         &&\qquad\qquad\qquad\qquad\qquad\quad \mbox{for some}\;k_1\in\N\;\mbox{and}\;k_2\in [1,n-1]\}. \nn\eea
Then the Betti numbers of the free loop space of $M$ are given by
\bea b_i(\Lm M)
&=& \rank H_i(\Lm M/S^1,\Lm^0 M/S^1;\Q)\nn\\
&=& \left\{\begin{array}{cccc}
    0, & \quad \mbox{if}\ i\ \mbox{is even or}\ i\le d-2,  \cr
    \left[\frac{i-(d-1)}{d}\right]+1, & \quad \mbox{if}\ i\in 2\N-1\;\mbox{and}\;d-1\le i < d-1+(n-1)d, \cr
    n+1, & \quad \mbox{if}\ i\in \Om(d,n), \cr
    n, & \quad \mbox{otherwise}. \end{array}\right.\lb{2.6}\eea
For any integer $k\ge dn-1$, we have
\be \sum_{i=0}^kb_i(\Lm M) = \frac{n(n+1)d}{2D}(k-(d-1)) - \frac{n(n-1)d}{4} + 1 + \Theta_{d,n}(k),  \lb{2.7}\ee
where
\bea &&\Theta_{d,n}(k) = \left\{\frac{D}{dn}\left\{\frac{k-(d-1)}{D}\right\}\right\}
          - \left(\frac{2}{d}+\frac{d-2}{dn}\right)\left\{\frac{k-(d-1)}{D}\right\}   \nn\\
&&\qquad\qquad\qquad - n\left\{\frac{D}{2}\left\{\frac{k-(d-1)}{D}\right\}\right\}
          - \left\{\frac{D}{d}\left\{\frac{k-(d-1)}{D}\right\}\right\}. \lb{2.8}\eea}

\subsection{Index iteration theory of closed geodesics}

In \cite{Lon99}, Y. Long established the basic normal form decomposition of symplectic matrices. Based on this result the precise iteration formulae of indices of symplectic paths has been established in \cite{Lon00}. Note that this index iteration formulae works for Morse indices
of iterated closed geodesics (cf. \cite{LL02}, \cite{Liu05} and Chapter 12 of \cite{Lon02}).
Since every closed geodesic on a compact simply-connected manifold $M$ is orientable, then by Theorem 1.1 of \cite{Liu05}, the Morse index of a closed geodesic $c$ on such manifold coincides with the index of a corresponding symplectic path.

As in \cite{Lon02}, denote by
\bea
N_1(\lm, a) &=& \left(
\begin{array}{ll}
\lm & a\cr
0 & \lm \\ \end{array}
\right), \qquad {\rm with\;}\lm=\pm 1, \; a\in\R, \lb{2.9}\\
H(b) &=&\left(
\begin{array}{ll}
b & 0\cr
0 & b^{-1}\\ \end{array}
\right), \qquad {\rm with\;}b\in\R\bs\{0, \pm 1\}, \lb{2.10}\\
R(\th) &=& \left(
  \begin{array}{ll}\cos\th & -\sin\th \cr
                           \sin\th & \cos\th\\ \end{array}
  \right), \qquad {\rm with\;}\th\in (0,\pi)\cup (\pi,2\pi), \lb{2.11}\\
N_2(e^{\th\sqrt{-1}}, B) &=&\left(
\begin{array}{ll}
 R(\th) & B \cr
 0 & R(\th)\\ \end{array}
 \right), \qquad {\rm with\;}\th\in (0,\pi)\cup (\pi,2\pi)\;\; {\rm and}\; \nn\\
        && \  B=\left(\begin{array}{ll} b_1 & b_2\cr
                                  b_3 & b_4\\ \end{array} \right)\; {\rm with}\; b_j\in\R, \;\;
                                         {\rm and}\;\; b_2\not= b_3. \lb{2.12}\eea
Here $N_2(e^{\th\sqrt{-1}}, B)$ is non-trivial if $(b_2-b_3)\sin\theta<0$, and trivial
if $(b_2-b_3)\sin\theta>0$ as defined in \cite{Lon02} and Definition 1.8.11 of \cite{Lon02}.

As in \cite{Lon02}, the $\diamond$-sum of any two real matrices is defined by
$$ \left(\begin{array}{ll} A_1 & B_1\cr C_1 & D_1\\ \end{array}\right)_{2i\times 2i}\diamond \left(\begin{array}{ll} A_2 & B_2\cr C_2 & D_2\\ \end{array}\right)_{2j\times 2j}
=\left(\begin{array}{llll} A_1 & 0 & B_1 & 0 \cr
                                   0 & A_2 & 0& B_2\cr
                                   C_1 & 0 & D_1 & 0 \cr
                                   0 & C_2 & 0 & D_2\\ \end{array}\right). $$

Let $J=\left(\begin{array}{ll}0&-I_n\cr
                 I_n&0\\ \end{array}\right)$ be the standard symplectic matrix, where $I_n$
is the identity matrix. As usual, the symplectic group is defined by $\Sp(2n) = \{M\in {\rm GL}(2n,\R)\,|\,M^TJM=J\}$.
For every $P\in\Sp(2n)$, the homotopy set $\Omega(P)$ of $P$ in $\Sp(2n)$ is defined by
$$ \Om(P)=\{Q\in\Sp(2n)\ |\ \sg(Q)\cap\U=\sg(P)\cap\U\equiv\Gamma\;\mbox{and}
                    \;\nu_{\om}(Q)=\nu_{\om}(P)\  \forall\ \om\in\Gamma\}, $$
where $\sg(P)$ denotes the spectrum of $P$,
$\nu_{\om}(P)\equiv\dim_{\C}\ker_{\C}(P-\om I)$ for $\om\in\U$.
The homotopy component $\Om^0(P)$ of $P$ in $\Sp(2n)$ is defined by
the path connected component of $\Om(P)$ containing $P$. Then the following decomposition theorem is proved in \cite{Lon99}
and \cite{Lon00}.

\medskip

{\bf Theorem 2.5.} (cf. Theorem 1.8.10, Lemma 2.3.5 and Theorem 8.3.1 of \cite{Lon02}) {\it  For every $P\in\Sp(2n)$, there
exists a continuous path $f:\,[0,1]\rightarrow\Om^0(P)$ such that $f(0)=P$ and
\bea f(1)
&=& N_1(1,1)^{\dm p_-}\,\dm\,I_{2p_0}\,\dm\,N_1(1,-1)^{\dm p_+}
  \dm\,N_1(-1,1)^{\dm q_-}\,\dm\,(-I_{2q_0})\,\dm\,N_1(-1,-1)^{\dm q_+} \nn\\
&&\dm\,N_2(e^{\aa_{1}\sqrt{-1}},A_{1})\,\dm\,\cdots\,\dm\,N_2(e^{\aa_{r_{\ast}}\sqrt{-1}},A_{r_{\ast}})
  \dm\,N_2(e^{\bb_{1}\sqrt{-1}},B_{1})\,\dm\,\cdots\,\dm\,N_2(e^{\bb_{r_{0}}\sqrt{-1}},B_{r_{0}})\nn\\
&&\dm\,R(\th_1)\,\dm\,\cdots\,\dm\,R(\th_r)\,\dm\,R(\th_{r+1})\,\dm\,\cdots\,\dm\,R(\th_{\bar{r}})\dm\,H(b)^{\dm h},\lb{2.13}\eea
where $\frac{\th_{j}}{2\pi}\in(0,1)\backslash\Q$ for $1\le j\le r$ and
$\frac{\th_{j}}{2\pi}\in((0,1)\setminus\{1/2\})\cap\Q$ for $r+1\le j\le \bar{r}$; $N_2(e^{\aa_{j}\sqrt{-1}},A_{j})$'s
are nontrivial and $N_2(e^{\bb_{j}\sqrt{-1}},B_{j})$'s are trivial, and non-negative integers
$p_\pm, p_0, q_\pm, q_0,\bar{r},r_\ast,r_0,h$ satisfy
\be p_- + p_0 + p_+ + q_- + q_0 + q_+ + \bar{r} + 2r_{\ast} + 2r_0 + h = n. \lb{2.14}\ee

Let $\ga\in\P_{\tau}(2n)=\{\ga\in C([0,\tau],\Sp(2n))\,|\,\ga(0)=I\}$, we extend
$\ga(t)$ to $t\in [0,m\tau]$ for every $m\in\N$ by
\be \ga^m(t)=\ga(t-j\tau)\ga(\tau)^j \qquad \forall\;j\tau\le t\le (j+1)\tau \;\;
               {\rm and}\;\;j=0, 1, \ldots, m-1. \lb{2.15}\ee

Denote the basic normal form decomposition of $P\equiv \ga(\tau)$ by (\ref{2.13}). Then we have
\bea i(\ga^m)
&=& m(i(\ga)+p_-+p_0-\bar{r}) + 2\sum_{j=1}^{\bar{r}}E\left(\frac{m\th_j}{2\pi}\right) - \bar{r}   \nn\\
&&  - p_- - p_0 - \frac{1+(-1)^m}{2}(q_0+q_+)
              + 2\sum_{j=1}^{r_{\ast}}\vf\left(\frac{m\aa_j}{2\pi}\right) - 2r_{\ast}. \lb{2.16}\eea}

Let
\be \M\equiv\{N_1(1,1); \;\;N_1(-1,a_2),\,a_2=\pm1;\;\;R(\th),
            \,\th\in[0,2\pi);\,H(-2)\}.               \lb{2.17}\ee
By Theorems 8.1.4-8.1.7 and 8.2.1-8.2.4 of \cite{Lon02}, we have

\medskip

{\bf Proposition 2.6.} {\it Every path $\ga\in\P_{\tau}(2)$ with end matrix homotopic
to some matrix in $\M$ has odd index $i(\ga)$. Every path $\xi\in\P_{\tau}(2)$
with end matrix homotopic to $N_1(1,-1)$ or $H(2)$, and every path $\eta\in\P_{\tau}(4)$
with end matrix homotopic to $N_2(\om,B)$ have even indices $i(\xi)$ and $i(\eta)$.}

\medskip

The common index jump theorem (cf. Theorem 4.3 of \cite{LZ02}) for symplectic paths has become one of the main tools in studying the multiplicity and stability of periodic orbits in Hamiltonian and symplectic dynamics. Recently, the following enhanced common index jump theorem has been obtained by Duan, Long and Wang in (\cite{DLW16a}).

\medskip

{\bf Theorem 2.7.} (cf. Theorem 3.5 of \cite{DLW16a}) {\it Let
$\gamma_k\in\mathcal{P}_{\tau_k}(2n)$ for $k=1,\cdots,q$ be a finite collection of symplectic paths.
Let $M_k=\ga_k(\tau_k)$. We extend $\ga_k$ to $[0,+\infty)$ by (\ref{2.15}) inductively. Suppose
\be  \hat{i}(\ga_k,1) > 0, \qquad \forall\ k=1,\cdots,q.  \lb{2.18}\ee
Then for any fixed integer $\bar{m}\in \N$, there exist infinitely many $(q+1)$-tuples
$(N, m_1,\cdots,m_q) \in \N^{q+1}$ such that for all $1\le k\le q$ and $1\le m\le \bar{m}$, there holds
\bea
\nu(\ga_k,2m_k-m) &=& \nu(\ga_k,2m_k+m) = \nu(\ga_k, m),   \lb{2.19}\\
i(\ga_k,2m_k+m) &=& 2N+i(\ga_k,m),                         \lb{2.20}\\
i(\ga_k,2m_k-m) &=& 2N-i(\ga_k,m)-2(S^+_{M_k}(1)+Q_k(m)),  \lb{2.21}\\
i(\ga_k, 2m_k)&=& 2N -(S^+_{M_k}(1)+C(M_k)-2\Delta_k),     \lb{2.22}\eea
where $S_{M_k}^\pm(\om)$ is the splitting number of $M_k$ at $\om$ (cf. Definition 9.1.4 of \cite{Lon02}) and
\be C(M_k)=\sum\limits_{0<\theta<2\pi}S^-_{M_k}(e^{\sqrt{-1}\theta}),\  \Delta_k = \sum_{0<\{m_k\th/\pi\}<\delta}S^-_{M_k}(e^{\sqrt{-1}\th}),\
 Q_k(m) = \sum_{e^{\sqrt{-1}\th}\in\sg(M_k),\atop\{\frac{m_k\th}{\pi}\} = \{\frac{m\th}{2\pi}\}=0}S^-_{M_k}(e^{\sqrt{-1}\th}). \lb{2.23}\ee
More precisely, by (4.10), (4.40) and (4.41) in \cite{LZ02} , we have
\bea m_k=\left(\left[\frac{N}{\bar{M}\hat i(\gamma_k, 1)}\right]+\chi_k\right)\bar{M},\quad 1\le k\le q,\lb{2.24}\eea
where $\chi_k=0$ or $1$ for $1\le k\le q$ and $\frac{\bar{M}\theta}{\pi}\in\Z$
whenever $e^{\sqrt{-1}\theta}\in\sigma(M_k)$ and $\frac{\theta}{\pi}\in\Q$
for some $1\le k\le q$.  Furthermore, for any fixed $M_0\in\N$, we may
further require  $M_0|N$, and for any $\epsilon>0$, we can choose $N$ and $\{\chi_k\}_{1\le k\le q}$ such that
\bea \left|\left\{\frac{N}{\bar{M}\hat i(\gamma_k, 1)}\right\}-\chi_k\right|<\epsilon,\quad 1\le k\le q.\lb{2.25}\eea}


\setcounter{equation}{0}
\section{Proof of Theorem 1.1}

Let $(M,F)$ be a compact simply-connected manifold with a bumpy, irreversible Finsler metric $F$ and satisfy $H^*(M;\Q)\cong T_{d,n+1}(x)$ of $\dim M=dn$, where integers $d\ge 2$ and $n\ge 1$. In order to prove Theorem 1.1, We make the following assumption

\medskip

{\bf (ECG)}\  {\it Suppose that all prime closed geodesics on $(M,F)$ are elliptic, and the total number of distinct closed geodesics is finite, denoted by $\{c_k\}_{k=1}^q$.}

\medskip

When the Finsler metric $F$ is bumpy, there holds $p_- + p_0 + p_+ + q_- + q_0 + q_+ + \bar{r}-r = 0 $ and there does not exist $\aa_j$ or $\bb_j$ which is the rational multiple of $\pi$ in (\ref{2.13}). Furthermore, if the closed geodesic is elliptic, then there also holds $h=0$. So by Theorem 2.5, under the assumption (ECG), for every prime closed geodesic $c_k$, $\forall\,1\le k\le q$, the basic normal form decomposition of the linearized Poincar\'{e} map $P_{c_k}$  possesses the following
form
\bea f_{c_k}(1)
&=& R(\th^k_{1})\,\dm\,\cdots\,\dm\,R(\th^k_{r_k}) \dm\,N_2(e^{\aa^k_{1}\sqrt{-1}},A^k_{1})\,\dm\,\cdots\,\dm\,N_2(e^{\aa^k_{r_{k\ast}}\sqrt{-1}},A^k_{r_{k\ast}})\nn\\   &&\dm\,N_2(e^{\bb^k_{1}\sqrt{-1}},B^k_{1})\,\dm\,\cdots\,\dm\,N_2(e^{\bb^k_{r_{k0}}\sqrt{-1}},B^k_{r_{k0}}), \nn\eea
where $\frac{\th^k_{j}}{2\pi}\in(0,1)\backslash\Q$ for $1\le j\le r_k$,
$\frac{\aa^k_{j}}{2\pi}\in (0,1)\backslash\Q$ for $1\le j\le r_{k\ast}$, $\frac{\bb^k_{j}}{2\pi}\in (0,1)\backslash\Q$ for $1\le j\le r_{k0}$,
and
\be r_k + 2r_{k\ast} + 2r_{k0} = dn - 1. \lb{3.1}\ee
Therefore, by (\ref{2.16}) we obtain the index iteration formula of $c_k^m$ for $1\le k\le q$
\be i(c_k^m) = m(i(c_k)-r_k) + 2\sum_{j=1}^{r_k}\left[\frac{m\th^k_j}{2\pi}\right] +r_k,\quad \nu(c_k^m)=0,
    \qquad\forall\ m\ge 1. \lb{3.2}\ee

\medskip

{\bf Claim 1.}\  {\it Under (ECG), for every $1\le k\le q$, there holds $i(c_k^m)=dn-1\ (\mod\,2)$, $\forall\ m\ge 1$.}

\medskip

In fact, by Proposition 2.6, (\ref{3.1}) and the homotopy invariant and the symplectic additivity of indices (cf. Theorem 6.2.7 of \cite{Lon02}), we obtain
\bea i(c_k)= dn-1 \quad (\mod\,2), \qquad 1\le k\le q. \lb{3.3}\eea

By (\ref{3.1}), (\ref{3.2}) and (\ref{3.3}), it yields
\bea i(c_k^{m+1})-i(c_k^m)&=&(i(c_k)-r_k) + 2\sum_{j=1}^{r_k}\left[\frac{(m+1)\th^k_j}{2\pi}\right]-2\sum_{j=1}^{r_k}\left[\frac{m\th^k_j}{2\pi}\right] \nn\\
&=&i(c_k)+2r_{k\ast} + 2r_{k0}-(dn-1)\quad (\mod\,2)\nn\\
&=& i(c_k)-(dn-1)\quad (\mod\,2)\nn\\
&=&0 \quad (\mod\,2),\qquad \forall\ 1\le k\le q,\quad m\ge 1. \lb{3.4}\eea
Now (\ref{3.3}) and (\ref{3.4}) finished the proof of Claim 1.

\medskip

Define the Morse-type numbers
\bea M_p=\sum_{k=1}^q M_p(k)\equiv\sum_{k=1}^q\#\{m\ge 1\ |\ i(c_k^m)=p, \;\ol{C}_p(E, c_k^m)\not= 0\},\quad p\in\Z,\lb{3.5}\eea
which, through the following Morse inequalities (cf. Theorem I.4.3 of \cite{Cha93}), relates the Betti numbers $b_p\equiv b_p(\Lm M)$ defined in Lemma 2.4
\bea M_p&\ge& b_p,\lb{3.6}\\
M_p - M_{p-1} + \cdots +(-1)^{p}M_0
&\ge& b_p - b_{p-1}+ \cdots + (-1)^{p}b_0, \quad\forall\ p\in\N_0.\lb{3.7}\eea

Since Claim 1 implies $i(c_k^m)-i(c_k)\in 2\Z,\ \forall\, m\ge 1$ and $1\le k\le q$.  By Proposition 2.1 we obtain
\bea M_p(k)=\#\{m\ge 1\ |\ i(c_k^m)=p\},\qquad \forall\ p\in\Z,\quad \forall\ 1\le k\le q, \lb{3.8}\eea
which, together with the fact $i(c_k^m)=dn-1\ (\mod\,2)$, $\forall\ m\ge 1$ from Claim 1, implies
\bea M_{p}=0, \qquad \forall\ p=dn\ (\mod 2),\ p\in\N_0.{\lb{3.9}}\eea

On the other hand, it follows from Lemma 2.4 that
\be b_{p}\equiv b_{p}(\Lm M)=0, \qquad \forall\ p=dn\ (\mod 2),\ p\in\N_0.\lb{3.10}\ee

Then by the Morse inequalities (\ref{3.6}) and (\ref{3.7}), we obtain
\be M_{p}=b_{p}, \qquad \forall\ p=dn-1\ (\mod 2),\ p\in\N_0.{\lb{3.11}}\ee

In summary, under the assumption (ECG), the Morse inequalities become the identities
\bea M_p = b_p,\quad \sum_{i=0}^p (-1)^i M_i=\sum_{i=0}^p (-1)^i b_i, \qquad \forall\ p\in\N_0.\lb{3.12}\eea

\medskip

{\bf Claim 2.} \ {\it Under (ECG), there holds $i(c_k)\ge d-1$ and the mean index $\hat{i}(c_k)>0$, $\forall\ 1\le k\le q$. Moreover, $\#\{1\le k\le q\ |\ i(c_k)=d-1\}=1$.}

\medskip

In fact, if there exists at least one closed geodesic $c_{k_0}$ such that $0\leq i(c_{k_0})< d-1$, then it follows from (\ref{3.5}) and (\ref{3.8}) that
$M_{i(c_{k_0})}=\sum_{k=1}^q M_{i(c_{k_0})}(k)\ge 1$. Note that $b_i=0, \forall\ 0\le i< d-1$ by Lemma 2.4. So by (\ref{3.12}) we get the following contradiction
\be 1\le M_{i(c_{k_0})} = b_{i(c_{k_0})} = 0.\lb{3.13}\ee
So there holds $i(c_k)\ge d-1$, $\forall\,1\le k\le q$.

Assume $\hat{i}(c_k)=0$ holds for some $1\le k\le q$, then it yields $i(c_k^m)=0$ for any $m\ge 1$ by Corollary 4.2 of \cite{LL02}. This contradicts to the fact $i(c_k)\ge d-1\ge 1$. Thus the mean index $\hat{i}(c_k)$ of every closed geodesic $c_k$ must be positive.

In addition, by Lemma 2.4 and (\ref{3.12}), it yields $M_{d-1}=b_{d-1}=1$. So by (\ref{3.8}), we have
\bea M_{d-1}=\sum_{k=1}^q M_{d-1}(k)=\sum_{k=1}^q\#\{m\ge 1\ |\ i(c_k^m)=d-1\}=1.\lb{3.14}\eea
So there exist only one pair $(k_0,m_0)$ with some $1\le k_0\le q$ and $m_0\in\N$ satisfying $i(c_{k_0}^{m_0})=d-1$. Note that there always holds $i(c_k^m)\ge i(c_k), \forall\,m\ge 1$ by the Bott-type formulae (cf. Theorem A of \cite{Bot56} and Theorem 9.2.1 of \cite{Lon02}). Since $i(c_{k_0})\ge d-1$, so the unique possibility is $m_0=1$, i.e., $\#\{1\le k\le q\ |\ i(c_k)=d-1\}=1$.

This completes the proof of Claim 2.

\medskip

By (\ref{3.2}) and Claim 2, it yields
\bea \hat{i}(c_k)=i(c_k)-r_k + \sum_{j=1}^{r_k}\frac{\th^k_j}{\pi}>0, \lb{3.15}\eea
which, together with (\ref{3.2}), implies that for any $m,l\in\N$, there holds
\bea  i(c_k^{m+l})-i(c_k^l)&=& m(i(c_k)-r_k)+ 2\sum_{j=1}^{r_k}\left[\frac{(m+l)\th^k_j}{2\pi}\right]-2\sum_{j=1}^{r_k}\left[\frac{l\th^k_j}{2\pi}\right] \nn\\
&\ge& m(i(c_k)-r_k)+\sum_{j=1}^{r_k}\frac{m\th^k_j}{\pi}-2r_k\nn\\
&=& m\hat{i}(c_k)-2r_k,\quad\forall\ 1\le k\le q.\lb{3.16}\eea

It follows from (\ref{3.15}) and (\ref{3.16}) that there exists sufficiently large $m\in\N$ such that $i(c_k^{m+l})-i(c_k^l)\ge 0, \ \forall\, l\ge 1$. So the positive
integer $\bar{m}$ defined by
\be \bar{m}=\max_{1\le k\le q}\left\{\min\{\hat{m}\in\N\ |\ i(c_k^{m+l})\ge i(c_k^l),\ \forall\ l\ge 1,\ m\ge\hat{m}\}\right\}
              \lb{3.17}\ee
is well-defined and finite.

For the integer $\bar{m}$ defined in (\ref{3.17}), it follows from Theorem 2.7
that there exist infinitely many $(q+1)$-tuples $(N, m_1, \ldots, m_q)\in\N^{q+1}$ such that
for any $1\le k\le q$, there holds
\bea
\bar{m}+2&\le&\min_{1\le k\le q}\{2m_k\},\lb{3.18}\\
i(c_k^{2m_k-m}) &=& 2N-i(c_k^m),\quad 1\le m\le\bar{m}, \lb{3.19}\\
i(c_k^{2m_k}) &=& 2N-C(M_k)+2\Delta_k,  \lb{3.20}\\
i(c_k^{2m_k+m}) &=& 2N+i(c_k^m),\quad 1\le m\le\bar{m}, \lb{3.21}\eea
where $M_k=P_{c_k}\in \Sp(2(dn-1))$ is the linearized Poincar\'e map of $c_k$. Note that in the bumpy case, $S^+_{M_k}(1)=0$ and $Q_k(m)=0$ holds for
all $m\in\N$.

\medskip

Next we continue the proof of Theorem 1.1 in two cases according to the parity of $d$.

\medskip

{\bf Case 1.} \ {\it $H^*(M;\Q)\cong T_{d,n+1}(x)$ with even $d\ge 2$.}

\medskip

We use three steps to carry out the proof of Theorem 1.1 in this case.

\medskip

{\bf Step 1}. {\it The existence of $\frac{dn(n+1)}{4}$ prime closed geodesics.}

\medskip

On one hand, since $i(c_k^m)\ge i(c_k)$ for any $m\ge 1$, so by (\ref{3.19}), (\ref{3.21}) and the
fact $i(c_k)\ge d-1\ge 1$ for $1\le k\le q$ from Claim 2, it yields
\bea
i(c_k^{2m_k-m}) &=& 2N-i(c_k^m)\le 2N-i(c_k)\le 2N-1,\quad 1\le m\le\bar{m}, \lb{3.22}\\
i(c_k^{2m_k+m}) &=& 2N+i(c_k^m)\ge 2N+i(c_k)\ge 2N+1,\quad 1\le m\le\bar{m}. \lb{3.23}\eea

On the other hand, for any $\bar{m}+1\le m<2m_k$, by the definition (\ref{3.17}) of $\bar{m}$ (note that $m-1\ge\bar{m}$), (\ref{3.19}), (\ref{3.21}) and $i(c_k)\ge 1$, we have
\bea i(c_k^{2m_k-m})&\le & i(c_k^{(2m_k-m)+(m-1)})=i(c_k^{2m_k-1})=2N-i(c_k)\le 2N-1, \lb{3.24}\\
i(c_k^{2m_k+m})&=&i(c_k^{(2m_k+1)+(m-1)})\ge i(c_k^{2m_k+1}) = 2N+i(c_k)\ge 2N+1,\lb{3.25}
\eea
where notice that the second inequality actually holds for any $m\ge\bar{m}+1$.

Therefore, for every $1\le k\le q$, by (\ref{3.22})-(\ref{3.25}), we obtained the following estimates
\bea
i(c_k^{2m_k-m}) &\le& 2N-i(c_k) \le 2N-1,\qquad \forall\ 1\le m < 2m_k, \lb{3.26}\\
i(c_k^{2m_k}) &=& 2N-C(M_k)+2\Delta_k, \lb{3.27}\\
i(c_k^{2m_k+m}) &\ge& 2N+i(c_k) \ge 2N+1,\qquad \forall\ m\ge 1. \lb{3.28}\eea

{\bf Claim 3.} {\it For $N\in\N$ in Theorem 2.7 satisfying (\ref{3.26})-(\ref{3.28}) and $2NB(d,n)\in2\N$,
we have
\be \sum_{1\le k\le q} 2m_k\gamma_{c_k}=2NB(d,n). \lb{3.29}\ee}

In fact, choose $\ep < \frac{1}{1+2\bar{M}\sum\limits_{1\le k\le q}|\ga_{c_k}|}$, by Theorem 2.3, (\ref{2.24}) and (\ref{2.25}) of Theorem 2.7, it yields
\bea \left|2NB(d,n)-\sum_{k=1}^q 2m_k\ga_{c_k}\right|
&=& \left|\sum_{k=1}^q\frac{2N\ga_{c_k}}{\hat{i}(c_k)}-\sum_{k=1}^q 2\ga_{c_k}
                    \left(\left[\frac{N}{\bar{M}\hat{i}(c_k)}\right]+\chi_k\right)\bar{M}\right| \nn\\
&\le& 2\bar{M}\sum_{k=1}^q |\ga_{c_k}|\left|\left\{\frac{N}{\bar{M}\hat{i}(c_k)}\right\}-\chi_k\right| \nn\\
&<& 2\bar{M}\ep\sum_{k=1}^q|\ga_{c_k}| \nn\\
&<& 1. \lb{3.30}\eea
Since every $2m_k\ga_{c_k}$ is an integer, Claim 3 is proved.

\medskip

Now by Proposition 2.1, Definition 2.2 and Claim 1, it yields
\bea \sum_{m=1}^{2m_k} (-1)^{i(c_k^m)} \dim \ol{C}_{i(c_k^{m})}(E,c_k^m)
&=& \sum_{i=0}^{m_k-1} \sum_{m=2i+1}^{2i+2} (-1)^{i(c_k^m)} \dim \ol{C}_{i(c_k^{m})}(E,c_k^m) \nn\\
&=& \sum_{i=0}^{m_k-1} \sum_{m=1}^{2} (-1)^{i(c_k^m)} \dim \ol{C}_{i(c_k^{m})}(E,c_k^m) \nn\\
&=& m_k \sum_{m=1}^{2} (-1)^{i(c_k^m)} \dim \ol{C}_{i(c_k^{m})}(E,c_k^m) \nn\\
&=& 2m_k\ga_{c_k},\qquad \forall\ 1\le k\le q, \lb{3.31}\eea
where the second equality follows from Proposition 2.1 and the fact $i(c_k^{m+2})-i(c_k^m)\in 2\Z$ for all
$m\ge 1$ from Claim 1, and the last equality follows from Proposition 2.1 and Definition 2.2.

By (\ref{3.28}) and Proposition 2.1, it follows that all $c_k^{2m_k+m}$'s with $m\ge1$ and $1\le k\le q$
have no contribution to the alternating sum $\sum_{p=0}^{2N}(-1)^p M_p$. Similarly again by Proposition 2.1
and (\ref{3.26}), all $c_k^{2m_k-m}$'s with $1\le m<2m_k$ and $1\le k\le q$ only have contribution
to $\sum_{p=0}^{2N}(-1)^p M_p$.

Thus for the Morse-type numbers $M_p$'s defined by (\ref{3.5}), by (\ref{3.31}) we have
\bea \sum_{p=0}^{2N}(-1)^p M_p
&=& \sum_{k=1}^{q}\ \sum_{1\le m\le 2m_k \atop i(c_k^{m})\le 2N} (-1)^{i(c_k^m)} \dim \ol{C}_{i(c_k^{m})}(E,c_k^m) \nn\\
&=& \sum_{k=1}^{q}\ \sum_{m=1}^{2m_k} (-1)^{i(c_k^m)}\dim\ol{C}_{i(c_k^{m})}(E,c_k^m) \nn\\
& & -\sum_{1\le k\le q \atop i(c_k^{2m_k})\ge 2N+1} (-1)^{i(c_k^{2m_k})} \dim \ol{C}_{i(c_k^{2m_k})}(E,c_k^{2m_k}) \nn\\
&=& \sum_{k=1}^{q} 2m_k\ga_{c_k}
  -\sum_{1\le k\le q \atop i(c_k^{2m_k})\ge 2N+1} (-1)^{i(c_k^{2m_k})} \dim \ol{C}_{i(c_k^{2m_k})}(E,c_k^{2m_k}).
      \lb{3.32}\eea

In order to precisely know the contribution of the iterate $c_k^{2m_k}$ of $c_k$ to the alternating sum $\sum_{p=0}^{2N}(-1)^p M_p(k)$ for $1\le k\le q$, we define
\bea
\mathfrak{N}_+ &=& \#\{1\le k\le q\ |\ i(c_k^{2m_k})\ge 2N+1\}, \lb{3.33}\\
\mathfrak{N}_- &=& \#\{1\le k\le q\ |\ i(c_k^{2m_k})\le 2N-1\}. \lb{3.34}\eea

Note that by Claim 1, there holds $i(c_k^{2m_k})-i(c_k)\in 2\N_0,\ i(c_k)\in 2\N-1$, $\forall\ 1\le k\le q$. Thus by Theorem 2.3, Claim 3, (\ref{3.12}), (\ref{3.32}), the definition of $\mathfrak{N}_+$ and Lemma 2.4, we have
\bea -\frac{Ndn(n+1)}{D}+\mathfrak{N}_+
&=& 2NB(d,n)+\mathfrak{N}_+  \nn\\
&=& \sum_{k=1}^q 2m_k\gamma_{c_k} + \mathfrak{N}_+ \nn\\
&=& \sum_{p=0}^{2N}(-1)^p M_p = \sum_{p=0}^{2N}(-1)^p b_p\nn\\
&=&-\sum_{2k-1=1}^{2N-1}b_{2k-1}(\Lm M) \nn\\
&=& -\frac{dn(n+1)}{2D}(2N-d)+\frac{dn(n-1)}{4}-1-\Theta_{d,n}(2N-1), \lb{3.35}\eea
which implies
\be \mathfrak{N}_+=\frac{d^2n(n+1)}{2D}+\frac{dn(n-1)}{4}-1-\Theta_{d,n}(2N-1).\lb{3.36}\ee

Note that we can assume that $N$ is a multiple of $D=d(n+1)-2$ by Theorem 2.7. So there holds $\{\frac{2N-1-(d-1)}{D}\}=1-\frac{d}{D}=\frac{dn-2}{D}$.
Then by (\ref{2.8}) we have
\bea \Theta_{d,n}(2N-1)
&=& \frac{dn-2}{dn}-\frac{2n+d-2}{dn}\left(1-\frac{d}{D}\right)
               -n\left\{\frac{dn-2}{2}\right\}-\left\{\frac{dn-2}{d}\right\} \nn\\
&=& \frac{(d-2)D+2d-d^2}{Dd}-\left\{-\frac{2}{d}\right\} = 1-\frac{2}{d}-\left\{-\frac{2}{d}\right\}-\frac{d-2}{D} \nn\\
&=& -\frac{d-2}{D}, \nn\eea
which, together with (\ref{3.36}), yields
\be \mathfrak{N}_+= \frac{d^2n(n+1)}{2D}+\frac{dn(n-1)}{4}-1+\frac{d-2}{D} = \frac{dn(n+1)}{4}. \lb{3.37}\ee

\medskip

{\bf Step 2.} {\it The existence of other $\frac{dn(n+1)}{4}$ prime closed geodesics.}

\medskip

Now we can apply Theorem 2.7 again and find another $(q+1)$-tuple $(N', m_1', \ldots, m_q')\in\N^{q+1}$ such
that similarly to (\ref{3.26})-(\ref{3.28}) for every $1\le k\le q$, there holds
\bea
i(c_k^{2m_k'-m})&\le& 2N'-i(c_k) \le 2N'-1,\qquad \forall\ 1\le m<2m_k, \lb{3.38}\\
i(c_k^{2m_k'})&=&2N'-C(M_k)+2\Delta'_k,  \lb{3.39}\\
i(c_k^{2m_k'+m})&\ge& 2N'+i(c_k) \ge 2N'+1,\qquad \forall\ m\ge 1,  \lb{3.40}\eea
where $$\Delta'_k = \sum_{0<\{m'_k\th/\pi\}<\delta}S^-_{M_k}(e^{\sqrt{-1}\th}), $$
and $\Delta_k$ and $\Delta'_k$ satisfy the following relationship (also cf. (42) in Theorem 2.8 of \cite{HW22})
\be  \Delta_k + \Delta'_k = C(M_k),\qquad 1\le k\le q. \lb{3.41}\ee

Similarly, note that $i(c_k^{2m_k'})-i(c_k)\in 2\N_0,\ i(c_k)\in 2\N-1,\ \forall\ 1\le k\le q$, we define
\bea
\mathfrak{N}'_+
&=& \#\{1\le k\le q\ |\ i(c_k^{2m_k'})\ge 2N'+1\}, \lb{3.42}\\
\mathfrak{N}'_-
&=& \#\{1\le k\le q\ |\ i(c_k^{2m_k'})\le 2N'-1\}. \lb{3.43}\eea

So by (\ref{3.39}) and (\ref{3.41}) it yields
\be i(c_k^{2m_k'}) = 2N'-C(M_k)+2(C(M_k)-\Delta_k)=2N'+C(M_k)-2\Delta_k. \lb{3.44}\ee

Note that when $i(c_k^{2m_k'})\ge 2N'+1$, it yields $C(M_k)-2\Delta_k\ge 1$ by (\ref{3.44}), then it follows from (\ref{3.27}) that $i(c_k^{2m_k}) = 2N-C(M_k)+2\Delta_k\le 2N-1$, and vice versa. So, by (\ref{3.27}), (\ref{3.39}), (\ref{3.44}) and the definitions of $\mathfrak{N}_\pm$ and $\mathfrak{N}'_\pm$, it yields
\be  \mathfrak{N}_{\pm} = \mathfrak{N}'_{\mp}. \lb{3.45}\ee

By Proposition 2.1 and (\ref{3.40}), it follows that all $c_k^{2m_k'+m}$'s with $m\ge1$ and $1\le k\le q$ have
no contribution to the alternating sum $\sum_{p=0}^{2N'}(-1)^p M_p$. Similarly also by Proposition 2.1 and
(\ref{3.38}), all $c_k^{2m_k'-m}$'s with $m\ge1$ and $1\le k\le q$ only have contribution to
$\sum_{p=0}^{2N'}(-1)^p M_p$.

Thus, through carrying out arguments similar to (\ref{3.35})-(\ref{3.37}), by Claim 3, the
definition of $\mathfrak{N}'_+$ and Lemma 2.4, together with (\ref{3.45}), we obtain
\be \mathfrak{N}_- = \mathfrak{N}'_+ = \frac{dn(n+1)}{4}. \lb{3.46}\ee

Since $i(c_k^{2m_k})\neq 2N$, $\forall\, 1\le k\le q$ by Claim 1, it yields $q=\mathfrak{N}_++\mathfrak{N}_-$ by the definitions (\ref{3.33}) and (\ref{3.34}) of $\mathfrak{N}_\pm$. Then by (\ref{3.37}) and (\ref{3.46}), we get
\bea q = \mathfrak{N}_++\mathfrak{N}_- = \frac{dn(n+1)}{4}+\frac{dn(n+1)}{4} = \frac{dn(n+1)}{2}.\lb{3.47}\eea
This completes the proof of Theorem 1.1 in this case.

\medskip

{\bf Case 2.} {\it $H^*(M;\Q)\cong T_{d,n+1}(x)$ with odd $d\ge 3$.}

\medskip

Note that in this case there holds $n=1$ and $d\ge 3$, which implies that $M$ is rationally homotopic to the sphere $S^d$.

Now we use two steps to carry out the proof of Theorem 1.1 in this case.

\medskip

{\bf Step 3.} {\it The existence of $(d-1)$ prime closed geodesics.}

\medskip

Firstly, note that there holds $i(c_k)\ge d-1\ge 2,\ \forall\, 1\le k\le q$ by Claim 2. Therefore, by Theorem 2.7 and some similar arguments as those above (\ref{3.26})-(\ref{3.28}), for $1\le k\le q$, we obtain
\bea
i(c_k^{2m_k-m}) &\le& 2N-i(c_k) \le 2N-2,\qquad \forall\ 1\le m < 2m_k, \lb{3.48}\\
i(c_k^{2m_k}) &=& 2N-C(M_k)+2\Delta_k, \lb{3.49}\\
i(c_k^{2m_k+m}) &\ge& 2N+i(c_k) \ge 2N+2,\qquad \forall\ m\ge 1. \lb{3.50}\eea

Therefore, similarly to the equation (\ref{3.32}), we have
\be \sum_{p=0}^{2N+1}(-1)^p M_p = \sum_{k=1}^{q} 2m_k\ga_{c_k}
  -\sum_{1\le k\le q \atop i(c_k^{2m_k})\ge 2N+2} (-1)^{i(c_k^{2m_k})}
          \dim\ol{C}_{i(c_k^{2m_k})}(E,c_k^{2m_k}). \lb{3.51}\ee

Define
\bea
\hat{\mathfrak{N}}_+ &=& \#\{1\le k\le q\ |\ i(c_k^{2m_k})\ge 2N+2\}, \lb{3.52}\\
\hat{\mathfrak{N}}_- &=& \#\{1\le k\le q\ |\ i(c_k^{2m_k})\le 2N-2\}. \lb{3.53}\eea

Note that by Claim 1, there holds $i(c_k^{2m_k})-i(c_k)\in 2\N_0,\ i(c_k)\in 2\N$, $\forall\ 1\le k\le q$. Then, by Theorem 2.3, (\ref{3.12}), Claim 3, (\ref{3.51}), the above definition of $\hat{\mathfrak{N}}_+$, and Lemma 2.4,
we have
\bea \frac{N(d+1)}{d-1}-\hat{\mathfrak{N}}_+
&=& 2NB(d,1)-\hat{\mathfrak{N}}_+ \,=\, \sum_{k=1}^q 2m_k\gamma_{c_k}-\hat{\mathfrak{N}}_+  \nn\\
&=& \sum_{p=0}^{2N+1}(-1)^p M_p = \sum_{p=0}^{2N+1}(-1)^p b_p\nn\\
&=& \sum_{2k=0}^{2N} b_{2k}(\Lm M)\nn\\
&=& \frac{N(d+1)}{d-1}-\frac{d-1}{2},\lb{3.54}\eea
which yields
\be  \hat{\mathfrak{N}}_+ = \frac{d-1}{2}.\lb{3.55}\ee

By Theorem 2.7 again, there exists another $(q+1)$-tuple $(N', m_1', \ldots, m_q')\in\N^{q+1}$ such
that (\ref{3.48})-(\ref{3.50}) hold with $N$ and $m_k$s being replaced by $N'$ and
$m_k'$s respectively. Then, denote by $\hat{\mathfrak{N}}'_+$ and $\hat{\mathfrak{N}}'_-$ the numbers similarly
defined by (\ref{3.52})-(\ref{3.53}) instead of using $N'$ and $m_k'$s. According to their definitions, obviously these numbers satisfy the relationship
\be \hat{\mathfrak{N}}_{\pm}= \hat{\mathfrak{N}}'_{\mp}. \lb{3.56}\ee

Similarly to the same arguments from (\ref{3.46}) and using the relationship (\ref{3.56}), we obtain
\be \hat{\mathfrak{N}}_- = \hat{\mathfrak{N}}'_+=\frac{d-1}{2}. \lb{3.57}\ee

Therefore from (\ref{3.55}) and (\ref{3.57}) it follows
\be  \hat{\mathfrak{N}}_++\hat{\mathfrak{N}}_- = \frac{d-1}{2}+\frac{d-1}{2} = d-1. \lb{3.58}\ee

\medskip

{\bf Step 4.} {\it The existence of other two prime closed geodesics.}

\medskip

Denote by $\{c_k\}_{k=1}^{d-1}$ the $(d-1)$ prime closed geodesics found in Step 3.
For these closed geodesics, there holds $i(c_k^{2m_k})\neq 2N$ by  the definitions of $\hat{\mathfrak{N}}_+$ and $\hat{\mathfrak{N}}_-$, which, together
with (\ref{3.48}) and (\ref{3.50}), yields
\be  i(c_k^m)\neq 2N,\qquad\forall\  m\ge 1,\  1\le k\le d-1.\lb{3.59}\ee
Then by Proposition 2.1 it yields
\be M_{2N}(k)=\#\{m\ge 1\ |\ i(c_k^m)=2N, \;\ol{C}_{2N}(E, c_k^m)\not= 0\} = 0,\quad \forall\ 1\le k\le d-1. \lb{3.60}\ee

Since $N$ can be chosen to be the multiple of $d-1$ by Theorem 2.7, by Lemma 2.4 and (\ref{3.12}), it yields $M_{2N}=b_{2N}(\Lm M)=2$. Therefore, by Proposition 2.1, (\ref{3.48}), (\ref{3.50}) and (\ref{3.60}) it yields
\bea 2&=&M_{2N}=\sum_{k=d}^q M_{2N}(k)\nn\\
&=&\sum_{k=d}^q\#\{m\ge 1\ |\ i(c_k^m)=2N, \;\ol{C}_{2N}(E, c_k^m)\not= 0\}\nn\\
&=&\sum_{k=d}^q\#\{m=2m_k\ |\ i(c_k^m)=2N, \;\ol{C}_{2N}(E, c_k^m)\not= 0\},\lb{3.61}\eea
where the final equality follows from (\ref{3.48}) and (\ref{3.50}).

Now it follows from (\ref{3.61}) and Proposition 2.1 that there exist exactly two prime closed geodesic $c_d$ and $c_{d+1}$ such that
\bea i(c_k^{2m_k})=2N, \qquad k=d, d+1,\lb{3.62}\eea which implies that $c_d$ and $c_{d+1}$ are distinct from $\{c_k\}_{k=1}^{d-1}$ by (\ref{3.59}). This completes the proof of Step 4.

Now notice that $dn-1=d-1$ is even in this case. By Claim 1, it yields $i(c_k^{2m_k})\not\in\{2N-1, 2N+1\}$, $\forall\, 1\le k\le q$. Therefore it yields $q=\hat{\mathfrak{N}}_++\hat{\mathfrak{N}}_-+2$ by (\ref{3.62}), the definitions (\ref{3.52}) and (\ref{3.53}) of $\hat{\mathfrak{N}}_\pm$. Then by (\ref{3.58}) and (\ref{3.62}), we get
\bea q = \hat{\mathfrak{N}}_++\hat{\mathfrak{N}}_-+2 = (d-1)+2 = d+1.\lb{3.63}\eea
This completes the proof of Theorem 1.1 in this case.  \hfill\hb

\bigskip

{\bf Acknowledgement.} The authors sincerely thank the referee for her/his careful reading and valuable comments.

\bibliographystyle{abbrv}

\begin{thebibliography}{100}
\addtolength{\itemsep}{-0.2ex}
\bibitem[Ano74]{Ano74} D. V. Anosov, Geodesics in Finsler geometry. {\it Proc.
I.C.M.} (Vancouver, B.C. 1974), Vol. 2. 293-297 Montreal (1975)
(Russian), {\it Amer. Math. Soc. Transl.} 109 (1977) 81-85.
\bibitem[Ban93]{Ban93} V. Bangert,  On the existence of closed geodesics on two-spheres.
{\it Internat. J. Math.} 4 (1993) 1-10.
\bibitem[BL10]{BL10} V. Bangert and Y. Long,   The existence of two closed
geodesics on every Finsler 2-sphere.  {\it Math. Ann.} 346 (2010) 335-366.
\bibitem[Bot56]{Bot56} R. Bott,  On the iteration of closed geodesics and the Sturm
intersection theory.  {\it Comm. Pure Appl. Math.} 9 (1956) 171-206.
\bibitem[BuK21]{BuK21} K. Burns and V. S. Matveev, Open problems and questions about geodesics. {\it Ergod. Th. \& Dynam. Sys.} 41 (2021) 641-684.
\bibitem[Cha93]{Cha93} K. C. Chang, Infinite Dimensional Morse Theory and Multiple Solution Problems, Birkh\"{a}user, Boston, 1993.
\bibitem[CHP19]{CHP19} D. Cristofaro-Gardiner, M. Hutchings and D. Pomerleano, Torsion contact forms
in three dimensions have two or infinitely many Reeb orbits. {\it Geom. Topol.}  23 (2019)  3601-3645.
\bibitem[Dua15]{Dua15} H. Duan, Non-hyperbolic closed geodesics on positively curved Finsler spheres. {\it J. Funct. Anal.} 269 (2015), no. 11, 3645-3662.
\bibitem[Dua16]{Dua16} H. Duan, Two elliptic closed geodesics on positively curved Finsler spheres. {\it J. Diff. Equ.} 260 (2016) 8388-8402.
\bibitem[DL16]{DL16} H. Duan and H. Liu, Multiplicity of closed geodesics on Finsler spheres with irrationally elliptic closed geodesics. {\it Sci. China Math.} 59 (2016) 531-538.
\bibitem[DL22]{DL22} H. Duan and H. Liu, The non-contractibility of closed geodesics on Finsler $RP^n$. {\it Acta Math. Sinica, English Series}. 38 (2022), no. 1, 1-21.
\bibitem[DuL07]{DuL07} H. Duan and Y. Long, Multiple closed geodesics on
bumpy Finsler $n$-spheres. {\it J. Diff. Equ.} 233 (2007) 221-240.
\bibitem[DuL10]{DuL10} H. Duan and Y. Long, The index growth and multiplicity
of closed geodesics. {\it J. Funct. Anal.} 259 (2010) 1850-1913.
\bibitem[DLW16a]{DLW16a} H. Duan, Y. Long and W. Wang, The enhanced common index jump theorem for symplectic paths and non-hyperbolic closed
geodesics on Finsler manifolds. {\it Calc. Var. and PDEs} 31 (2016) 1-28.
\bibitem[DLW16b]{DLW16b} H. Duan, Y. Long and W. Wang, Two closed geodesics on compact simply connected bumpy Finsler manifolds. {\it J. Diff. Geom.} 104 (2016) 275-289.
\bibitem[DLZ20]{DLZ20} H. Duan, Y. Long and C. Zhu, Index iteration theories for periodic orbits: old and new. {\it Nonlinear Anal.} 201 (2020), 111999, 26 pp.
\bibitem[Fra92]{Fra92} J. Franks, Geodesics on $S\sp 2$ and periodic points of annulus homeomorphisms.
{\it Invent. Math.} 108 (1992) 403-418.
\bibitem[GrM69]{GrM69} D. Gromoll and W. Meyer, Periodic geodesics on compact Riemannian
manifolds. {\it J. Diff. Geom.}  3 (1969) 493-510.
\bibitem[GG20]{GG20} V. L. Ginzburg and B. Z. G\"{u}rel, Lusternik-Schnirelmann theory and closed Reeb orbits. {\it Math. Z.} 295 (2020) 515-582.
\bibitem[GGM18]{GGM18} V. L. Ginzburg, B. Z. G\"{u}rel and L. Macarini, Multiplicity of closed Reeb orbits on prequantization bundles. {\it Israel J. Math.} 228 (2018) 407-453.
\bibitem[GM21]{GM21} V. L. Ginzburg and L. Macarini, Dynamical convexity and closed orbits on symmetric spheres. {\it Duke Math. J.} 170 (2021) 1201-1250.
\bibitem[HW22]{HW22} M. Hamid and W. Wang, A symmetric property in the enhanced common index jump theorem with applications to the closed geodesic problem. {\it Discrete Contin. Dyn. Syst.} 42 (2022), 1933-1948.
\bibitem[Hin84]{Hin84} N. Hingston,  Equivariant Morse theory and closed geodesics. {\it J. Diff. Geom.} 19 (1984) 85-116.
\bibitem[Hin93]{Hin93} N. Hingston, On the growth of the number of closed geodesics on the two-sphere. {\it Inter. Math. Research Notices}. 9 (1993) 253-262.
\bibitem[HR13]{HR13} N. Hingston and H.-B. Rademacher, Resonance for loop homology of spheres. {\it J. Diff. Geom.} 93 (2013) 133-174.
\bibitem[HWZ98]{HWZ98} H. Hofer, K. Wysocki and E. Zehnder, The dynamics on three-dimensional strictly convex energy surfaces. {\it Ann. of Math.} 148 (1998) 197-289.
\bibitem[HWZ03]{HWZ03} H. Hofer, K. Wysocki and E. Zehnder, Finite energy foliations of tight three shperes and Hamiltonian dynamics. {\it Ann. of Math.} 157 (2003) 125-255.
\bibitem[Kat73]{Kat73} A. B. Katok, Ergodic properties of degenerate
integrable Hamiltonian systems. {\it Izv. Akad. Nauk SSSR.} 37 (1973) (Russian),
{\it Math. USSR-Izv.} 7 (1973) 535-571.
\bibitem[Kli78]{Kli78} W. Klingenberg, Lectures on Closed Geodesics. Springer.
Berlin. 1978.
\bibitem[Liu05]{Liu05} C. Liu,  The relation of the Morse index of closed geodesics
with the Maslov-type index of symplectic paths. {\it Acta Math. Sinica, English
Series.} 21 (2005) 237-248.
\bibitem[LL02]{LL02} C. Liu and Y. Long,  Iterated index formulae for closed
geodesics with applications. {\it Science in China.} 45 (2002) 9-28.
\bibitem[Liu17]{Liu17} H. Liu, The Fadell-Rabinowitz index and multiplicity of non-contractible closed geodesics on Finsler $RP^n$.
{\it J. Diff. Equ.} 262 (2017) 2540-2553.
\bibitem[LX17]{LX17} H. Liu and Y. M. Xiao, Resonance identity and multiplicity of non-contractible closed geodesics on Finsler
$RP^n$. {\it Adv. Math.} 318 (2017) 158-190.
\bibitem[LLX18]{LLX18} H. Liu, Y. Long and Y. Xiao, The existence of two non-contractible closed geodesics
on every bumpy Finsler compact space form. {\it Disc. Cont. Dyna. Syst.} 38 (2018) 3803-3829.
\bibitem[Lon99]{Lon99} Y. Long,  Bott formula of the Maslov-type index theory. {\it Pacific J. Math.} 187 (1999) 113-149.
\bibitem[Lon00]{Lon00} Y. Long,  Precise iteration formulae of the Maslov-type index theory and ellipticity of closed characteristics.
{\it Adv. Math.} 154 (2000) 76-131.
\bibitem[Lon02]{Lon02} Y. Long,  Index Theory for Symplectic Paths with Applications. Progress in Math. 207, Birkh\"auser. Basel. 2002.
\bibitem[Lon06]{Lon06} Y. Long,  Multiplicity and stability of closed geodesics on Finsler 2-spheres. {\it J. Eur. Math.
Soc.} 8 (2006) 341-353.
\bibitem[LD09]{LD09}  Y. Long and H. Duan, Multiple closed geodesics on 3-spheres. {\it Adv. Math.} 221 (2009) 1757-1803.
\bibitem[LW08]{LW08}  Y. Long and W. Wang, Stability of closed geodesics on Finsler 2-spheres. {\it J. Funct. Anal.} 255 (2008) 620-641.
\bibitem[LZ02]{LZ02} Y. Long and C. Zhu,  Closed characteristics on compact convex hypersurfaces in $\R^{2n}$.  {\it Ann. of Math.}
155 (2002) 317-368.
\bibitem[Rad89]{Rad89} H.-B. Rademacher,  On the average indices of closed
geodesics. {\it J. Diff. Geom.} 29 (1989) 65-83.
\bibitem[Rad92]{Rad92} H.-B. Rademacher, Morse Theorie und geschlossene
Geodatische. {\it Bonner Math. Schriften} Nr. 229 (1992).
\bibitem[Rad07]{Rad07} H.-B. Rademacher, Existence of closed geodesics on
positively curved Finsler manifolds. {\it Ergod. Th. \& Dynam. Sys.} 27 (2007) 957--969.
 \bibitem[Rad10]{Rad10} H.-B. Rademacher, The second closed geodesic on Finsler
spheres of dimension $n>2$. {\it Trans. Amer. Math. Soc. } 362 (2010) 1413-1421.
\bibitem[She01]{She01} Z. Shen, Lectures on Finsler Geometry. World Scientific.
Singapore. 2001.
\bibitem[Tai10]{Tai10} I. A. Taimanov, The type numbers of closed geodesics. {\it Regul. Chaotic Dyn.} 15 (2010) 84-100.
\bibitem[VS76]{ViS76} M. Vigu\'e-Poirrier and D. Sullivan,  The homology theory of the closed
geodesic problem. {\it J. Diff. Geom.} 11 (1976) 633-644.
\bibitem[Wan08]{Wan08} W. Wang, Closed geodesics on positively curved Finsler spheres.
{\it Adv. Math.} 218 (2008) 1566-1603.
\bibitem[Wan12]{Wan12} W. Wang, On a conjecture of Anosov. {\it Adv. Math.} 230 (2012) 1597-1617.
\bibitem[Wan13]{Wan13} W. Wang, On the average indices of closed geodesics on positively curved Finsler spheres.
{\it Math. Ann.} 355 (2013) 1049-1065.
\bibitem[Zil82]{Zil82} W. Ziller, Geometry of the Katok examples. {\it Ergod. Th. \& Dynam. Sys.} 3 (1982) 135-157.

\end{thebibliography}

\end{document}